\documentclass[11pt,reqno]{amsart}
\usepackage{amsthm,amssymb,amsfonts}
\usepackage[nosumlimits]{mathtools}
\usepackage{mathrsfs}
\usepackage{graphicx}
\usepackage{tikz-cd}
\usepackage{enumerate}
\usepackage{enumitem}
\usepackage{subcaption}

\usepackage[margin=1in]{geometry}

\usepackage{hyperref,xcolor}
\hypersetup{
    colorlinks,
    linkcolor={red!50!black},
    citecolor={blue!50!black},
    urlcolor={blue!80!black}
}
\usepackage{adjustbox}

\let\oldlhd\unlhd
\renewcommand*{\unlhd}{\mathrel{\mskip0.1mu \oldlhd \mskip0.1mu}}

\theoremstyle{plain}
\newtheorem{theorem}{Theorem}[section]
\newtheorem{proposition}[theorem]{Proposition}
\newtheorem{lemma}[theorem]{Lemma}

\theoremstyle{definition}

\newcommand{\F}{\mathbb F}
\newcommand{\E}{\mathbb E}

\DeclareMathOperator{\codim}{codim}

\DeclareMathOperator{\Hom}{Hom}

\DeclareMathOperator{\bias}{bias}

\DeclareMathOperator{\Prob}{Pr}
\DeclareMathOperator{\arank}{AR}
\DeclareMathOperator{\prank}{PR}
\DeclareMathOperator{\rk}{rank}
\DeclareMathOperator{\im}{im}

\begin{document}
\title{partition rank and analytic rank over finite field}
\date{}
\author[Q-Y.~Chen]{Qiyuan~Chen}
\address{Korteweg-de Vries Institute for Mathematics, University of Amsterdam,
P.O. Box 94248, 1090 GE Amsterdam, The Netherlands}
\email{q.chen3@uva.nl}
\author[K.~Ye]{Ke Ye}
\address{State Key Laboratory of Mathematical Sciences, Academy of Mathematics and Systems Science, Chinese Academy of Sciences, Beijing 100190, China}
\email{keyk@amss.ac.cn}
\begin{abstract}
We prove that the partition rank and the analytic rank of a tensor over any finite field are equivalent up to a constant depending only on the order of the tensor. The proof combines Mili\'cevi\'c's theorem on dense multilinear varieties with the differential Schur complement of Moshkovitz and Zhu. The new ingredient is a method for approximately eliminating multilinear equations.
\end{abstract}
\maketitle
\section{Introduction}\label{sec:intro}
\subsection{Background and main results}
A basic question about polynomials over finite fields is how the bias of a polynomial constrains its algebraic structure. For multilinear
forms, this question can be expressed precisely in terms of two notions of rank.
Let $d\geq 2$, let $\F=\mathbb F_q$, let $V_1,\ldots,V_d$ be finite-dimensional
$\F$-vector spaces, and let $T\colon V_1\times\cdots\times V_d\to\F$
be a multilinear form. Its \emph{partition rank} $\prank_{\F}(T)$ is
the least number of summands in an expression
\[
 T(x_1,\ldots,x_d)=\sum_{\nu=1}^s
 f_\nu(x_{I_\nu})g_\nu(x_{[d]\setminus I_\nu}),
 \qquad \varnothing\ne I_\nu\subsetneq[d],
\]
where $[d]=\{1,\ldots,d\}$, $x_I=(x_i)_{i\in I}$, and both factors
are multilinear in the indicated variables. Naslund \cite{Naslund} introduced
this notion in the context of the polynomial method. The
\emph{analytic rank}, introduced by Gowers and Wolf \cite{GowersWolf}, is
\[
 \arank_{\F}(T)=-\log_q\bias_{\F}(T),\qquad
 \bias_{\F}(T)=\E_{x_1,\ldots,x_d}\chi(T(x_1,\ldots,x_d)),
\]
where $\chi$ is any nontrivial additive character of $\F$. For a
multilinear form, this expectation is positive and independent of
$\chi$. Thus partition rank measures the length of a shortest decomposition,
whereas analytic rank measures the quasirandomness of the tensor. Both ranks agree with matrix rank when $d=2$.

The inequality $\arank_{\F}(T)\leq\prank_{\F}(T)$ follows immediately from
the subadditivity of analytic rank proved in \cite{lovett2018analytic}. The reverse comparison
has been the subject of a long sequence of results. The work of Green
and Tao \cite{GreenTao}, together with refinements by Kaufman and Lovett
\cite{KaufmanLovett} and by Bhowmick and Lovett \cite{BhowmickLovett},
gives a qualitative bound for partition rank in terms of analytic rank.
Mili\'cevi\'c \cite{milicevic2019polynomial} and Janzer \cite{janzer2019polynomial} independently obtained polynomial bounds. For trilinear
forms, Adiprasito, Kazhdan, and Ziegler \cite{adiprasito2021schmidt} and, independently, Cohen and Moshkovitz \cite{CMBilinear} proved a linear bound. Both proofs use Derksen's stability result for slice rank \cite{derksen2022g}.
Lampert \cite{Lampert2025} later gave an elementary proof of this bound without using geometric invariant theory.
For forms of arbitrary order, Cohen and Moshkovitz \cite{CM} proved a linear
bound when the field size is sufficiently large in terms of the
analytic rank. Moshkovitz and Zhu \cite[Theorem~2]{MZLocal} sharpened the dependence on the field size and obtained
the bound
\begin{equation}\label{eq:previous-bound}
 \prank_{\F}(T)\ll_d
 \arank_{\F}(T)\bigl(1+\log_q(\arank_{\F}(T)+1)\bigr)
\end{equation}
over every finite field. Here and below, the implicit constant in
$\ll_d$ depends only on $d$.

The partition rank versus analytic rank conjecture
\cite{lovett2018analytic,KZApplications,adiprasito2021schmidt} asserts that the logarithmic factor in \eqref{eq:previous-bound} can be removed. In this paper, we prove this conjecture.

\begin{theorem}\label{thm:main}
For every integer $d\geq2$, there is a constant $C_d$ such that, for
every finite field $\F$ and every $d$-linear form $T$ over $\F$,
\[
 \arank_{\F}(T)\leq\prank_{\F}(T)
 \leq C_d\,\arank_{\F}(T).
\]
\end{theorem}

In particular, the field size need not grow with the analytic rank.
We do not attempt to optimize the dependence of $C_d$ on $d$.
In terms of the bias, Theorem~\ref{thm:main} reads
\[
 q^{-\prank_{\F}(T)}\leq\bias_{\F}(T)
 \leq q^{-\prank_{\F}(T)/C_d}.
\]
This gives an exponential estimate for equidistribution in terms of
partition rank. For example, let $T_1,\ldots,T_a$ be $d$-linear forms
on the same spaces, and suppose that $\prank_{\F}(\lambda_1T_1+\cdots+\lambda_aT_a)\geq R$
for every $0\ne\lambda\in\F^a$. Fourier inversion and
Theorem~\ref{thm:main} give, for every $c\in\F^a$,
\[
 \left|\Prob_x\bigl[(T_1(x),\ldots,T_a(x))=c\bigr]-q^{-a}\right|
 \leq q^{-R/C_d}.
\]
Indeed, each nontrivial Fourier coefficient is the bias of such a
linear combination. This is a direct application to counting rational
solutions of systems of multilinear equations.

The comparison also enters higher-order Fourier analysis through
polarization. The standard reductions in \cite[\S7.3]{CM} and
\cite[\S\S8.2, 8.5]{MZLocal}, applied to Theorem~\ref{thm:main},
give a polynomial bound in the degree-$d$ versus degree-$(d-1)$ case of the Gowers inverse theorem
for polynomial phases over $\mathbb F_p$ with $p>d$. Specifically, a degree-$d$
polynomial phase with $U^d$-norm at least $\delta$ has correlation at least
$\delta^{O_d(1)}$ with some polynomial phase of degree at most $d-1$. Bias and rank estimates also play a role in the
geometry and universality of high-rank polynomial maps
\cite{KZApplications} and in the study of Reed--Muller codes
\cite{BhowmickLovett}. These connections explain the importance of
bounds that are independent of the number of variables.

Our proof of Theorem~\ref{thm:main} uses the partial decomposition strategy proposed by Moshkovitz and Zhu \cite[\S7]{MZSchur}. We remove a form whose partition rank is bounded linearly in the original analytic rank and
reduce the analytic rank of the remainder by a fixed factor. The
remainder is a form on the original spaces, so the procedure can be
iterated. More precisely, we prove the following theorem.

\begin{theorem}[Partial decomposition]\label{thm:partial}
For each $d\geq3$, there are constants $Q_d,L_d,B_d$ depending only on
$d$ with the following property. Let $|\F|=q\geq Q_d$ and let $T$
be a $d$-linear form with $r=\arank_{\F}(T)$. There exists a
$d$-linear form $T'$ on the same spaces such that, for $D=2^{d-3}$,
\begin{equation}\label{eq:partial}
 \prank_{\F}(T-T')\leq L_d r,\qquad
 \arank_{\F}(T')\leq
 \min\{r/2,\ B_d(r+1)q^{-1/D}\}.
\end{equation}
\end{theorem}

\subsection{Overview of the proof}\label{subsec:overview}
We first explain how Theorem~\ref{thm:partial} gives
Theorem~\ref{thm:main}. When $d=2$, both ranks equal the matrix rank. For $d\geq3$, suppose $q\geq Q_d$, increasing $Q_d$ if necessary so
that every nonzero $d$-linear form has analytic rank at least $1/2$
(Lemma~\ref{lem:bias}). Starting with $T_0=T$, apply
Theorem~\ref{thm:partial} successively to obtain remainders $T_i$ with
\[
 \arank_{\F}(T_i)\leq2^{-i}\arank_{\F}(T),\qquad
 \prank_{\F}(T_i-T_{i+1})\leq L_d\arank_{\F}(T_i).
\]
The remainders eventually vanish. Summing the bounds gives $\prank_{\F}(T)\leq2L_d\arank_{\F}(T)$. For smaller finite fields, the same conclusion up to a constant depending only on $d$ follows by the standard scalar extension argument of \cite[Lemmas~8.1 and~8.2]{MZLocal}. Thus it suffices to prove
Theorem~\ref{thm:partial}.

Write $k=d-2$. For fixed $x$, the last two blocks of $T$ define a bilinear form with matrix $M_T(x)$, so that
\[
 T(x,z,w)=\langle M_T(x)z,w\rangle,\qquad
 x\in X=X_1\times\cdots\times X_k.
\]
The proof combines three ingredients.

\smallskip\noindent
\emph{Dense multilinear varieties.}
The first ingredient, Proposition~\ref{prop:fibers}, follows from Mili\'cevi\'c's theorem on dense multilinear varieties \cite{Milicevic}. In the form used here, this theorem gives a subvariety defined by $O_d(r+1)$ multilinear equations inside any multilinear variety of density $q^{-r}$. Applied to the equation $M_T(x)z=0$, it gives a subvariety of $X$ defined by $O_d(r+1)$ multilinear equations on which $M_T$ has rank $O_d(r+1)$.

\smallskip\noindent
\emph{A relative differential Schur complement.}
The second ingredient, Proposition~\ref{prop:relativeschur}, adapts the proof of \cite[Theorem~1.1]{MZSchur}. Suppose that $V\subseteq X$ is defined by $s$ multilinear equations and that $e=\max_{x\in V}\rk_{\F} M_T(x)$. By refining the system to $O_d(s)$ equations and subtracting a form of partition rank $O_d(e)$ from $T$, we obtain a remainder whose maximum matrix rank on the refined variety is $O_d(e/(q-k))$ when $q>k$.

\smallskip\noindent
\emph{Approximate elimination of multilinear equations.}
The third ingredient, Proposition~\ref{prop:all}, uses a bound on the matrix rank over $V$ to obtain a bound on the analytic rank over the whole space. With $s$ and $e$ as above, we find a form $S$ with $\prank_{\F}(T-S)\ll_d s$. When $k=1$, we have $\arank_{\F}(S)\leq e$. When $k\geq2$, we have $\arank_{\F}(S)\ll_d s/m+m^{2^{k-1}-1}e$ for every $m$ in the range of Proposition~\ref{prop:all}\textup{(b)}.
\subsection{Organization}
Section~\ref{sec:notation} collects the notation and preliminary
facts. Section~\ref{sec:input} contains the density result and the
relative differential Schur complement. Section~\ref{sec:eliminate}
develops the approximate elimination argument.
Section~\ref{sec:main} combines these ingredients to prove
Theorem~\ref{thm:partial}.
\section{Notation and preliminaries}\label{sec:notation}
Throughout this paper, $\F=\mathbb F_q$ is a finite field, and all vector spaces are finite-dimensional over $\F$ unless otherwise specified. Write $[h]=\{1,\ldots,h\}$. Given spaces $X_1,\ldots,X_h$ and $I\subseteq[h]$, set $X_I=\prod_{i\in I}X_i$ and write $x_I=(x_i)_{i\in I}$. In a partial evaluation, the entries retain their original positions. We write $A\ll_h B$ if $A\leq c_h B$ for a constant $c_h$ depending only on $h$. All expectations and probabilities are taken with respect to the uniform distribution on the relevant spaces.
\subsection{Partition rank and analytic rank}
We regard an order-$d$ tensor as a multilinear form $T\colon V_{[d]}\to\F$. For $d\geq 2$, its \emph{partition rank} $\prank_{\F}(T)$ is the least integer $s$ such that
\[
 T(x)=\sum_{\nu=1}^s f_\nu(x_{I_\nu})\,g_\nu(x_{[d]\setminus I_\nu}),
 \qquad \varnothing\ne I_\nu\subsetneq[d],
\]
where both factors are multilinear in the indicated blocks. Concatenating decompositions gives the subadditivity of $\prank_{\F}$. For a nontrivial additive character $\chi$ of $\F$, the \emph{bias} and the \emph{analytic rank} of $T$ are given by
\[
 \bias_{\F}(T)=\E_{x\in V_{[d]}}\chi(T(x)),\qquad
 \arank_{\F}(T)=-\log_q\bias_{\F}(T).
\]

\begin{lemma}\label{lem:bias}
For $d\geq2$, the bias of $T$ is positive. For any $i\in[d]$, it equals
the probability that the induced linear functional on $V_i$ is zero
when all other blocks are chosen independently and uniformly \cite[(1)]{lovett2018analytic}.
Analytic rank is subadditive \cite[Theorem~1.5]{lovett2018analytic} and does not
increase under linear pullbacks in the individual blocks
\cite[Lemma~1.4.7]{BDR}. Moreover, by \cite[Theorem~1.7(i)]{lovett2018analytic},
\begin{equation}\label{eq:easy}
 \arank_{\F}(T)\leq\prank_{\F}(T).
\end{equation}
For bilinear forms, partition rank equals matrix rank by rank
factorization, and analytic rank has the same value
\cite[Example~1.4]{lovett2018analytic}. If $T\ne0$ and $q\geq(d-1)^2$, then
$\arank_{\F}(T)\geq1/2$.
\end{lemma}

For the last claim, choose $v_i\in V_i$ with $T(v_1,\ldots,v_d)\ne0$.
On the lines $\F v_i$, the restricted tensor has a single nonzero entry.
By monotonicity and \cite[Theorem~1.7(ii)]{lovett2018analytic} applied to this
singleton, $\arank_{\F}(T)\geq1-\log_q(d-1)\geq1/2$.

\subsection{Multilinear equations and multilinear approximation}\label{subsec:mixed}
On $X=X_{[h]}$, a multilinear equation with support $I$, where
$\varnothing\ne I\subseteq[h]$, is an equation $\gamma(x_I)=0$ for
a scalar multilinear form $\gamma$ on $X_I$. For a finite system
$\Gamma$ of such equations, we denote its set of common $\F$-rational
zeros by $V(\Gamma)$ and call it a \emph{multilinear variety}. We write
$|\Gamma|$ for the number of equations. Every nonempty fiber obtained
by fixing all but one block is a linear subspace. Indeed, the equations
involving the free block become homogeneous linear equations in that
block, and the others are already satisfied.

For $a,x\in X$ and $t=(t_1,\ldots,t_h)$, write $a+t\odot x=(a_i+t_i x_i)_{i=1}^h$.
Following \cite[\S2.1]{MZSchur}, we write $f[a_S]$ for the function
obtained by fixing the blocks in $S$ at $a_S$. We work with rational functions
in the function field $\F(X)$. Such a function is regular at
$a$ if it can be written as a quotient whose denominator is nonzero at $a$.
Moshkovitz and Zhu \cite[Definition~3.1]{MZSchur} define the
\emph{multilinear approximation} of a rational function $f$ of multidegree
$(1,\ldots,1)$ that is regular at $a$ by
\[
 [f]_a(x)=
 \left.\frac{\partial}{\partial t_1}\cdots
       \frac{\partial}{\partial t_h}\,f(a+t\odot x)\right|_{t=0}.
\]
We use the same formula for any rational function regular at $a$ and
apply it entrywise to matrices. The derivatives are formal. Since each
$t_i$ is differentiated only once, the value is the coefficient of
$t_1\cdots t_h$ in the expansion at $t=0$, in every characteristic.

For fixed $a$ and $J\subseteq[h]$, write $J^c=[h]\setminus J$ and set
\[
 f_J(x_J)=\bigl[f[a_{J^c}]\bigr]_{a_J}(x_J).
\]
Thus $f_J$ is obtained by fixing the blocks in $J^c$ and taking the
multilinear approximation in the blocks in $J$. In particular,
$f_{[h]}=[f]_a$ and $f_\varnothing=f(a)$. If $f$ is already
multilinear, then $f_J(x_J)=f[a_{J^c}](x_J)$.
For a rational matrix $R(u_1,\ldots,u_h)$ in scalar variables that is regular at the origin, we
write $[R]_0(\mathbf 1)$ for this approximation at the origin evaluated
at the all-ones vector. Its value is the coefficient of $u_1\cdots u_h$.
When substituting indeterminates into a blockwise affine map, we use
its polynomial expression of total degree at most one in each block.
We write $\rk_{\F(t)}$ for rank over the rational function field $\F(t)$.

\begin{lemma}\label{lem:coefficient rule}
Let $A$ and $B$ be rational matrices regular at $a\in X$. For every
$J\subseteq[h]$, $A_J(x_J)$ is multilinear in $x_J$. This applies
in particular to blockwise affine matrices, sums and products of
compatible matrices, and $A^{-1}$ when $A(a)$ is invertible.
For a product, we have
\[
 (AB)_J(x_J)=\sum_{I\subseteq J} A_I(x_I)\,B_{J\setminus I}(x_{J\setminus I}).
\]
\end{lemma}
\begin{proof}
Expand each entry of $A(a+v)$ as a formal power series in the block
coordinates of $v$. After the substitution $v_i=t_ix_i$, the mixed
derivative with respect to the variables $t_i$ for $i\in J$, evaluated at $t=0$, selects precisely
the terms of degree one in each of those blocks and degree zero in the others.
These terms are multilinear in $x_J$. Sums and products of regular matrices
are regular, as is $A^{-1}=\operatorname{adj}(A)/\det A$ when
$\det A(a)\ne0$. In the Leibniz rule, the summand indexed by $I$ occurs when
each derivative with index in $I$ acts on $A$ and each remaining
derivative acts on $B$.
\end{proof}

We use the partial evaluations from the proof of
\cite[Theorem~1.2, \S3.2]{ChenYe} in the following finite-field form.
\begin{lemma}\label{lem:refine}
Let $\Gamma$ be a system of multilinear equations on $X_{[h]}$, and let
$a\in V(\Gamma)$. There is a system $\Gamma_a$ of at most
$(2^h-1)|\Gamma|$ multilinear equations with $V(\Gamma_a)\subseteq
V(\Gamma)$ such that
\[
 \prod_{i=1}^h \operatorname{span}_{\F}\{a_i,x_i\}\subseteq V(\Gamma)
 \qquad(x\in V(\Gamma_a)).
\]
\end{lemma}
\begin{proof}
For each $\gamma\in\Gamma$ with support $I$ and each nonempty $J\subseteq I$, include the equation
$\gamma[a_{I\setminus J}](x_J)=0$ in $\Gamma_a$.
This gives at most $2^{|I|}-1$ equations per $\gamma$ and includes
$\gamma(x_I)=0$ when $J=I$. If $x\in V(\Gamma_a)$, all these mixed
evaluations vanish. The evaluation for $J=\varnothing$ also vanishes
since $a\in V(\Gamma)$. By multilinearity,
$\gamma((s_i a_i+t_i x_i)_{i\in I})$ is a linear combination of these
evaluations for every choice of scalars $s_i,t_i\in\F$. This proves the inclusion.
\end{proof}

For $k=d-2\geq1$, we denote the last two blocks of a form by $Z$ and $W$ and set $X=\prod_{i=1}^kX_i$. The form determines a multilinear map
\[
 M_T\colon X\longrightarrow\Hom_{\F}(Z,W^*),\qquad
 \langle M_T(x)z,w\rangle=T(x,z,w),
\]
where $\langle\lambda,w\rangle=\lambda(w)$. We write $M=M_T$ when $T$ is fixed.
In particular, $\rk_{\F}M(x)$ is the ordinary rank of the bilinear form obtained by fixing $x$.
\section{Dense subvarieties and the relative differential Schur complement}\label{sec:input}
\subsection{Dense multilinear varieties and matrix fibers}\label{sec:density}
Mili\'cevi\'c states the following theorem for prime fields. However, each step
of his proof in \cite[\S\S2--3]{Milicevic} applies over an arbitrary
finite field $\F_q$. The directional convolutions use only addition,
and the density of a subspace of codimension $r$ over $\F_q$ is $q^{-r}$.

\begin{theorem}[{\cite[Theorem~2]{Milicevic}}]\label{thm:weak}
For each $h\geq 1$ there is a constant $K_h$ with the following property.
Let $G_1,\ldots,G_h$ be finite-dimensional vector spaces over $\F$,
where $|\F|=q$, and let $V\subseteq G_1\times\cdots\times G_h$ be a
multilinear variety of density $q^{-\rho}$. There is a set $\Gamma$ of
at most $K_h(\rho+1)$ multilinear equations such that $V(\Gamma)\subseteq V$.
\end{theorem}
The following proposition is a consequence of Theorem~\ref{thm:weak}.
\begin{proposition}\label{prop:fibers}
Let $X=X_1\times\cdots\times X_k$, where $k\geq1$, and let
$T\colon X\times Z\times W\to\F$ be multilinear. Write
$T(x,z,w)=\langle M(x)z,w\rangle$ and $r=\arank_{\F}(T)$.
There is a set $\Gamma$ of $s$ multilinear equations on $X$ such that
\[
 s+\max_{x\in V(\Gamma)}\rk_{\F}M(x)\leq K_{k+1}(r+1).
\]
\end{proposition}

\begin{proof}
Let $P=\{(x,z)\in X\times Z:M(x)z=0\}$. By the definition of analytic rank, the density of $P$ in $X\times Z$ is $q^{-r}$.
Theorem~\ref{thm:weak} therefore gives a multilinear subvariety $U\subseteq P$ defined
by $N\leq K_{k+1}(r+1)$ equations. Let $\Gamma$ consist of the
equations defining $U$ whose supports omit $Z$, and write $s=|\Gamma|$.
The remaining $N-s$ equations have the form $N(x)z=0$, where $N(x)$ is
an $(N-s)\times\dim Z$ matrix whose rows are multilinear on their respective supports in $X$.
For every $x_0\in V(\Gamma)$, we have
\[
 \ker N(x_0)\subseteq\ker M(x_0),
\]
since $U\subseteq P$. It follows that
\[
 \rk_{\F}M(x_0)\leq \rk_{\F}N(x_0)\leq N-s,
\]
where the last inequality holds because $N(x_0)$ has $N-s$ rows.
Taking the maximum over $x_0\in V(\Gamma)$ gives
\[
 s+\max_{x_0\in V(\Gamma)}\rk_{\F}M(x_0)
 \leq N\leq K_{k+1}(r+1).
\]
\end{proof}
\subsection{A relative differential Schur complement}\label{sec:schur}
The differential Schur complement introduced by Moshkovitz and Zhu
\cite[Definition~4.3]{MZSchur} is the multilinear approximation of a rational matrix. They used it to prove the following result.
\begin{theorem}[{\cite[Theorem~1.1]{MZSchur}}]\label{thm:PRvsCR}
Let $M\colon X_1\times\cdots\times X_k\to\Hom_{\F}(Y,Z)$ be a multilinear matrix. Then
\[
 \prank(M)\ll_k\max_{x\in X_1\times\cdots\times X_k}\rk_{\F}M(x),
\]
where $\prank(M)$ denotes the matrix partition rank as defined in \cite{MZSchur}.
\end{theorem}
We prove a relative version of Theorem~\ref{thm:PRvsCR} in Proposition~\ref{prop:relativeschur} by modifying the proof in \cite{MZSchur}. The following lemma is based on the proof of \cite[Lemma~5.2]{MZSchur}.
\begin{lemma}\label{lem:line}
Let $I\subseteq\F$ have size $n$, and let $C(t)$ be a polynomial matrix whose entries have degree at most $h<n$. If $\rk_{\F} C(c)\leq e$ for all $c\in I$, then
\[
 \rk_{\F(t)} C(t)-e\leq\frac{he}{n-h}.
\]
\end{lemma}
\begin{proof}
Put $g=\rk_{\F(t)}C(t)$. We may assume that $g>e$. Choose a $g\times g$ submatrix $M$ whose determinant $f(t)$ is nonzero. Then $\deg f\leq hg$. For each $c\in I$, put $r_c=\rk_{\F}M(c)\leq e$. By constant elementary row operations, we may arrange that $g-r_c$ rows of $M$ vanish at $c$. Every entry in these rows is divisible by $t-c$, so $(t-c)^{g-r_c}$ divides $f(t)$.
Counting roots with multiplicity gives $n(g-e)\le\sum_{c\in I}(g-r_c)\leq\deg f\leq hg$.
This gives the bound.
\end{proof}
The following lemma adapts the argument of \cite[Proposition~3.7]{MZSchur}. We include the proof for completeness.
\begin{lemma}\label{lem:boolean}
Let $f(u)$ be a rational matrix in $h\ge1$ variables, regular at the origin. For $S\subseteq[h]$, let $\mathbf 1_S$ be its indicator vector. Then
\[
 \rk_{\F}\bigl([f]_0(\mathbf 1)\bigr)
 \leq\sum_{S\subseteq[h]}\rk_{\F(t)}f(t\mathbf 1_S).
\]
If $f(0)=0$, the term for $S=\varnothing$ is zero.
\end{lemma}
\begin{proof}
Consider the alternating sum along lines from \cite[Lemma~3.5]{MZSchur},
\[
 F(t)=\sum_{S\subseteq[h]}(-1)^{h-|S|}f(t\mathbf 1_S).
\]
We claim that $F(t)=[f]_0(\mathbf 1)t^h+O(t^{h+1})$. Let $u^{\alpha}$ be a monomial in the Taylor expansion of an entry of $f$. If $u^{\alpha}$ omits $u_i$, then $u^{\alpha}(t\mathbf 1_S)=u^{\alpha}(t\mathbf 1_{S\cup\{i\}})$, so
\[
 \sum_{S\subseteq[h]}(-1)^{h-|S|}u^{\alpha}(t\mathbf 1_S)
 =\sum_{S\subseteq[h]\setminus\{i\}}(-1)^{h-|S|}
 \bigl(u^{\alpha}(t\mathbf 1_S)-u^{\alpha}(t\mathbf 1_{S\cup\{i\}})\bigr)=0.
\]
Thus $u^{\alpha}$ can contribute to $F(t)$ only if $\operatorname{supp}(\alpha)=[h]$. The only such monomial of degree at most $h$ is $u_1\cdots u_h$, whose coefficient matrix is $[f]_0(\mathbf 1)$. This proves the claim.

It follows that $t^{-h}F(t)$ is regular at $0$ and has value $[f]_0(\mathbf 1)$ there. Therefore
\[
 \rk_{\F}([f]_0(\mathbf 1))\leq\rk_{\F(t)}F(t)
 \leq\sum_{S\subseteq[h]}\rk_{\F(t)}f(t\mathbf 1_S).
\]
\end{proof}
Lemma~\ref{lem:boolean} gives the following analogue of \cite[Proposition~4.5]{MZSchur}.
\begin{lemma}\label{lem:schur-coefficient}
Let $0\in I\subseteq\F$ with $m=|I|>k\geq1$. Let $L(u)$ be a matrix
of polynomials of degree at most one in each of $u_1,\ldots,u_k$, and write
\[
 L(u)=\begin{pmatrix}A&B\\ C&D\end{pmatrix},\qquad b=\rk_{\F}L(0).
\]
Assume the $b\times b$ block $A(0)$ is invertible and
$\rk_{\F}L(c\mathbf 1_J)\leq b$ for $c\in I$ and
$\varnothing\ne J\subseteq[k]$.  Put $g_J=\rk_{\F(t)}L(t\mathbf 1_J)$ for each
nonempty $J$. Then the Schur complement $E(u)=D-CA^{-1}B$ satisfies
\begin{equation}\label{eq:schur-rank-sum}
 \rk_{\F}\bigl([E]_0(\mathbf 1)\bigr)
 \leq\sum_{\varnothing\ne J\subseteq[k]}(g_J-b)
 \leq\sum_{\varnothing\ne J\subseteq[k]}\frac{|J|b}{m-|J|}
 \leq\frac{k2^{k-1}b}{m-k}.
\end{equation}
\end{lemma}
\begin{proof}
The invertibility of $A(0)$ ensures that $E$ is regular at zero. Block elimination at zero gives
\[
 \rk_{\F}E(0)=\rk_{\F}L(0)-b=0.
\] Since $A(t\mathbf 1_J)$ is invertible over $\F(t)$, block
elimination of $L(t \mathbf 1_J)$ implies
$\rk_{\F(t)}E(t\mathbf 1_J)=g_J-b$. 

Applying
Lemma~\ref{lem:boolean} to $E$ yields
\[
 \rk_{\F}\bigl([E]_0(\mathbf 1)\bigr)
 \leq \sum_{\varnothing\ne J\subseteq[k]}\rk_{\F(t)}E(t\mathbf 1_J)
 =\sum_{\varnothing\ne J\subseteq[k]}(g_J-b).
\]
The term for the empty subset vanishes because $E(0)=0$.
For nonempty $J$, the entries of $L(t\mathbf 1_J)$ have degree at most
$|J|$. By hypothesis,
$\rk_{\F}L(c\mathbf 1_J)\leq b$ at all $m$ values $c\in I$.
By Lemma~\ref{lem:line},
\[
 \rk_{\F(t)} L(t\mathbf 1_J)-b\leq\frac{|J|b}{m-|J|}.
\]
Summing these bounds gives the second inequality in \eqref{eq:schur-rank-sum}.
The last inequality follows from $m-|J|\geq m-k$ and $\sum_{J\subseteq[k]}|J|=k2^{k-1}$,
since each index occurs in $2^{k-1}$ subsets.
\end{proof}
The following proposition gives a relative version of the differential Schur complement estimate. By adding multilinear equations and subtracting a tensor of bounded partition rank, we reduce the maximum matrix rank on the resulting variety.
\begin{proposition}\label{prop:relativeschur}
Let \(X=\prod_{i=1}^{k}X_i\), where \(1\le k<q\), and let
\(T(x,z,w)=\langle M(x)z,w\rangle\) be a multilinear form on
\(X\times Z\times W\). For a set \(\Gamma\) of \(s\) multilinear
equations on $X$, let
\[
 b=\max_{x\in V(\Gamma)}\rk_{\F}M(x).
\]
There exist a set \(\Gamma'\) of multilinear equations on $X$ and a
multilinear form \(S(x,z,w)=\langle M_S(x)z,w\rangle\) on
\(X\times Z\times W\) such that
\[
 V(\Gamma')\subseteq V(\Gamma),\qquad
 |\Gamma'|\le (2^k-1)s,\qquad \prank_{\F}(T-S)\le 2^k b.
\]
Moreover, for every \(x\in V(\Gamma')\),
\begin{equation}\label{eq:schur-rank}
 \rk_{\F}M_S(x)
 \le \sum_{\varnothing\ne J\subseteq[k]}\frac{|J|b}{q-|J|}
 \le \frac{k2^{k-1}b}{q-k}.
\end{equation}
\end{proposition}
\begin{proof}
If $b=0$, take $S=T$ and $\Gamma'=\Gamma$. Otherwise, choose
$a\in V(\Gamma)$ with $\rk_{\F} M(a)=b$. By Lemma~\ref{lem:refine},
the refinement $\Gamma'$ at $a$ satisfies
$|\Gamma'|\leq(2^k-1)s$ and $V(\Gamma')\subseteq V(\Gamma)$. Moreover, \begin{equation}\label{eq:schur-lines}
 a+u\odot x\in V(\Gamma)\qquad(x\in V(\Gamma'),\ u\in\F^k).
\end{equation}

After choosing suitable bases in $Z$ and $W$, we may assume that
\[
 M=\begin{pmatrix}A&B\\ C&D\end{pmatrix},
 \qquad M(a)=\begin{pmatrix}I_b&0\\0&0\end{pmatrix}.
\]
Here $A$ is the $b\times b$ block. Write $E=D-CA^{-1}B$. Following \cite[Definition~4.3]{MZSchur}, set
\[
 M/A=\begin{pmatrix}0&0\\0&E\end{pmatrix},\qquad
 M_S=[M/A]_a,\qquad S(x,z,w)=\langle M_S(x)z,w\rangle.
\]
The entries of $M/A$ have multidegree $(1,\ldots,1)$, since those
of $A^{-1}$ have multidegree $(-1,\ldots,-1)$.
Thus $M_S$ is the differential Schur complement and is multilinear in the blocks of $X$.

By \cite[Proposition~4.4]{MZSchur}, $\prank_{\F}(T-S)\le 2^k b$. It suffices to show \eqref{eq:schur-rank}.
Fix any $x\in V(\Gamma')$. Set $L(u)=M(a+u\odot x)$.
The Schur complement of $L(u)$ with respect to $A(a+u\odot x)$ is
$E_x(u)=E(a+u\odot x)$, so $[E_x]_0(\mathbf 1)=[E]_a(x)$.
By \eqref{eq:schur-lines}, $\rk_{\F}L(u)=\rk_{\F}M(a+u\odot x)\leq b$ for every $u\in\F^k$. Applying Lemma~\ref{lem:schur-coefficient} to $L$ with $I=\F$, we obtain from \eqref{eq:schur-rank-sum} that
\[
 \rk_{\F} M_{S}(x)=\rk_{\F} [E_x]_0(\mathbf 1)
 \leq\frac{k2^{k-1}b}{q-k}.
\]
\end{proof}
\section{Elimination of multilinear equations}\label{sec:eliminate}
\subsection{Eliminating a system of linear equations defined by a polynomial matrix}\label{subsec:kernel}

Let $A(u)\colon Y\to E$ be a polynomial matrix of degree at most one in
each of $u_1,\ldots,u_\ell$, where $\ell\geq1$. Put $e=\rk A(0)$.
Choosing suitable bases of $Y$ and $E$, we can write
\[
 A(u)=\begin{pmatrix}Q&R\\ S&T\end{pmatrix},\qquad
 A(0)=\begin{pmatrix}I_e&0\\0&0\end{pmatrix}.
\]
For $y=(v,k)$, the equation $A(u)y=0$ gives $Qv+Rk=0$.
Since $Q(0)=I_e$, we can solve this equation for $v$ near zero while keeping $k$
fixed. The resulting projection is
\begin{equation}\label{eq:kernel-projector}
 P(u)(v,k)=(-Q^{-1}Rk,k),\qquad
 P(u)=\begin{pmatrix}0&-Q^{-1}R\\0&I\end{pmatrix}.
\end{equation}
Substituting $v=-Q^{-1}Rk$ into $Sv+Tk=0$ gives $(T-SQ^{-1}R)k=0$. Thus
\[
 AP=\begin{pmatrix}0&0\\ 0&T-SQ^{-1}R\end{pmatrix}
\]
is the Schur complement $A/Q$.
For $J\subseteq[\ell]$,
we write $P_J(\mathbf 1)$ for the mixed derivative at zero in the variables indexed by $J$, with $P_\varnothing(\mathbf 1)=P(0)$.

In this section, we write
$u^J=\prod_{j\in J}u_j$, with $u^\varnothing=1$.

\begin{lemma}[Interpolation on kernels]\label{lem:kernel-interpolation}
Let $A$ and $P$ be the polynomial and rational matrices defined above, respectively. Let $0\in I\subseteq\F$ with $m=|I|>\ell$, and set $\Omega=\{c\mathbf 1_D:c\in I,\ D\subseteq[\ell]\}$. Assume that $\rk A(\xi)\leq e=\rk A(0)$ for every $\xi\in\Omega$. Let $G$ consist of the vectors $y\in Y$ for which there exists a tuple
$(v_\xi)_{\xi\in\Omega}$, with $v_\xi\in\ker A(\xi)$, satisfying
\begin{equation}\label{eq:kernel-moments}
 P_{J^c}(\mathbf 1)y=\sum_{\xi\in\Omega}\xi^J v_\xi
 \qquad\text{for every }J\subseteq[\ell].
\end{equation}
Then $G$ is a linear subspace of $Y$ and
\begin{equation}\label{eq:kernel-codim}
 \codim_Y G\leq
 \sum_{\varnothing\ne D\subseteq[\ell]}\frac{|D|e}{m-|D|}
 \leq\frac{\ell 2^{\ell-1}e}{m-\ell}.
\end{equation}
\end{lemma}

\begin{proof}
First we show that $G$ is a linear subspace of $Y$.
Write $U=\bigoplus_{J\subseteq[\ell]}Y$ and define linear maps
\begin{align*}
 L\colon\bigoplus_{\xi\in\Omega}\ker A(\xi)&\longrightarrow U,
 & (v_{\xi})_{\xi}&\longmapsto \biggl(\sum_{\xi\in\Omega}\xi^{J}v_{\xi}\biggr)_{J},\\
 C\colon Y&\longrightarrow U,
 & y&\longmapsto \bigl(P_{J^c}(\mathbf 1)y\bigr)_{J}.
\end{align*}
By definition, $G=C^{-1}(\im L)$, so $G$ is a linear subspace of $Y$.

Let $\F[u_1,\ldots,u_\ell]_{\mathrm{aff}}$ denote the space of polynomials of degree at most one in each variable. We identify $U^*$ with $Y^*\otimes_{\F}\F[u_1,\ldots,u_\ell]_{\mathrm{aff}}$ via the map $\phi\colon(\lambda_J)_J\mapsto\sum_{J\subseteq[\ell]}\lambda_J\otimes u^J$. Under this identification, each element of $U^*$ is a polynomial row vector.

For each $\lambda=(\lambda_J)_J\in U^*$, we have
\[
 \lambda\bigl(L((v_{\xi})_\xi)\bigr)
 =\sum_{\xi\in\Omega}\sum_{J\subseteq[\ell]}\xi^{J}\lambda_{J}(v_{\xi})
 =\sum_{\xi\in\Omega}\phi(\lambda)(\xi)v_{\xi}.
\]
Thus $\lambda\in(\im L)^\perp$ if and only if $\phi(\lambda)(\xi)v=0$ for every $\xi\in\Omega$ and every $v\in\ker A(\xi)$. Equivalently, $\phi(\lambda)(\xi)$ lies in the row space of $A(\xi)$ for every $\xi\in\Omega$.

Now $Cy\in\im L$ if and only if $\lambda(Cy)=0$ for every $\lambda\in(\im L)^\perp$. We have
\[
 \lambda(Cy)
 =\sum_{J\subseteq[\ell]}\lambda_{J}\bigl(P_{J^c}(\mathbf 1)y\bigr)
 =\sum_{J\subseteq[\ell]}\phi(\lambda)_{J}(\mathbf 1)P_{J^c}(\mathbf 1)y
 =[\phi(\lambda)P]_0(\mathbf 1)y,
\]
where the last equality follows from Lemma~\ref{lem:coefficient rule}.

Choose a basis of $(\im L)^\perp$ and let $\Lambda(u)$ be the $\codim_U(\im L)\times\dim Y$ matrix whose rows are the corresponding polynomial row vectors. Consequently,
\[
 G=\ker\bigl([\Lambda P]_0(\mathbf 1)\bigr),\qquad
 \codim_Y G=\rk_{\F}\bigl([\Lambda P]_0(\mathbf 1)\bigr).
\]
Write $\Lambda=(\Lambda_1,\Lambda_2)$ and set
\[
 M=\begin{pmatrix}Q&R\\ \Lambda_1&\Lambda_2\end{pmatrix}.
\]
For every $\xi\in\Omega$, each row of $M(\xi)$ lies in the row space
of $A(\xi)$. Hence
\[
 \rk_{\F}M(\xi)\leq \rk_{\F}A(\xi)\leq e.
\]
Since $Q(0)=I_e$, we also have $\rk_{\F}M(0)=e$. Moreover,
\[
 M/Q
 =\begin{pmatrix}0&0\\0&\Lambda_2-\Lambda_1Q^{-1}R\end{pmatrix}
 =\binom{0}{\Lambda P}.
\]
Applying Lemma~\ref{lem:schur-coefficient} to $M$ gives
\[
 \codim_Y G
 =\rk_{\F}\bigl([\Lambda P]_0(\mathbf 1)\bigr)
 \leq\sum_{\varnothing\ne D\subseteq[\ell]}\frac{|D|e}{m-|D|}
 \leq\frac{\ell 2^{\ell-1}e}{m-\ell}.
\]
\end{proof}

\begin{lemma}\label{lem:analytic quotient}
Let $T\colon X_1\times\cdots\times X_k\times Y\to\F$ be a multilinear form, and let $Y'$ be a subspace of $Y$ of codimension $s$. Let $T'$ be the restriction of $T$ to $X_1\times\cdots\times X_k\times Y'$ induced by the inclusion $Y'\hookrightarrow Y$. Then $\arank(T')\geq\arank(T)-s$.
\end{lemma}
\begin{proof}
We have $T=T'+\sum_{i=1}^s x_i\otimes y_i$ for some $x_i\in X_1^*\otimes\cdots\otimes X_k^*$ and $y_i\in Y^*$. The claim follows from Lemma~\ref{lem:bias}.
\end{proof}

The following proposition allows us to eliminate multilinear equations one layer at a time at a linear cost in partition rank.
\begin{proposition}\label{prop:kernel}
Let $X=X_1\times\cdots\times X_\ell$ and
$Z=Z_1\times\cdots\times Z_p$, where $\ell\geq1$ and $p\geq2$. Let $\Gamma$ be a system of multilinear equations on $X$, and let
$B\colon X\to\Hom_{\F}(Y,E)$ be affine in each block of $X$.
Let $H\colon X\times Y\times Z\to\F$ be multilinear. Suppose that
\[
 e=\max_{x\in V(\Gamma)}\rk_{\F}B(x),\qquad
 \arank_{\F}\bigl(H(x,y,\cdot)\bigr)\leq b
 \quad (x\in V(\Gamma),\ y\in\ker B(x)).
\]
There exist a system $\Gamma'$ of multilinear equations on $X$ and a multilinear form
$\widehat H$ on $X\times Y\times Z$ such that
\[
 V(\Gamma')\subseteq V(\Gamma),\qquad
 |\Gamma'|\leq(2^\ell-1)|\Gamma|,\qquad
 \prank_{\F}(H-\widehat H)\leq 2^\ell e.
\]
For every $x\in V(\Gamma')$ and every integer $m$ with $\ell<m\leq q$, we have
\begin{equation}\label{eq:kernel-tail}
 \arank_{\F}\bigl(\widehat H(x,\cdot,\cdot)\bigr)
 \leq \frac{\ell 2^{\ell-1}e}{m-\ell}
       +\bigl(1+(2^\ell-1)(m-1)\bigr)b .
\end{equation}
\end{proposition}

\begin{proof}
Choose $a\in V(\Gamma)$ with $\rk_{\F}B(a)=e$. In suitable bases, we have
\[
 B=\begin{pmatrix}Q&R\\S&T\end{pmatrix},\qquad
 B(a)=\begin{pmatrix}I_e&0\\0&0\end{pmatrix}.
\]
Define $P=\begin{pmatrix}0&-Q^{-1}R\\0&I\end{pmatrix}$ as above. Let $\mathcal H=\Hom_\F(Z_1\otimes\cdots\otimes Z_p,\F)$. We identify $H$ with a polynomial matrix $h\colon X\to\Hom_\F(Y,\mathcal H)$. Set $\widehat h=[hP]_a$, and let $\widehat H$ be the tensor associated with $\widehat h$. This is well defined since $P$ is regular at $a$. Write $h=(h_1,h_2)$. Since $h$ is multilinear in the blocks of $X$,
we have $[h]_a=h$. Hence
\[
 h-\widehat h=[h(I_Y-P)]_a.
\]
If $e=0$, then $P=I_Y$ and $\widehat h=h$. For $e>0$, set
\[
 D=\begin{pmatrix}I_e&Q^{-1}R\end{pmatrix}.
\]
The definition of $P$ gives
\[
 h(I_Y-P)=h_1D.
\]
For $1\leq i\leq e$, let $c_i$ be the $i$th column of $h_1$ and
let $d_i$ be the $i$th row of $D$. Thus $c_i$ takes values in
$\mathcal H$, $d_i$ takes values in $Y^*$, and both are regular at $a$.
By Lemma~\ref{lem:coefficient rule},
\[
 (h-\widehat h)(x)
 =\sum_{i=1}^e\sum_{J\subseteq[\ell]}
   (c_i)_J(x_J)\,(d_i)_{J^c}(x_{J^c}).
\]
After evaluation at $y\in Y$ and $z\in Z$, the summand indexed by
$(i,J)$ is
\[
 \bigl((c_i)_J(x_J)\bigr)(z)\,
 \bigl((d_i)_{J^c}(x_{J^c})\bigr)(y).
\]
The first factor is multilinear in $x_J,z_1,\ldots,z_p$, and the
second is multilinear in $x_{J^c},y$. These two sets of blocks form
a nontrivial partition for every $J\subseteq[\ell]$, including
$J=\varnothing$ and $J=[\ell]$. Thus
\[
 \prank_{\F}(H-\widehat H)\leq 2^\ell e.
\]

Next, choose $\Gamma'$ as in Lemma~\ref{lem:refine}. It suffices to prove \eqref{eq:kernel-tail}. Fix $x\in V(\Gamma')$ and write $x_\xi=a+\xi\odot x$. Let $G\subseteq Y$ be the subspace in Lemma~\ref{lem:kernel-interpolation} satisfying \eqref{eq:kernel-moments} and \eqref{eq:kernel-codim}. For each $y\in G$, there exist vectors $v_\xi\in\ker B(x_\xi)$ such that $P_{J^c}(x_{J^c})y=\sum_{\xi\in\Omega}\xi^Jv_\xi$ for every $J\subseteq[\ell]$. Hence, by Lemma~\ref{lem:coefficient rule},
\begin{align*}
 \widehat h(x)y
 &=\sum_{J\subseteq[\ell]} h_J(x_J)P_{J^c}(x_{J^c})y\\
 &=\sum_{\xi\in\Omega}\biggl(\sum_{J\subseteq[\ell]}\xi^J h_J(x_J)\biggr)v_\xi
   =\sum_{\xi\in\Omega} h(a+\xi\odot x)v_\xi,
\end{align*}
where the last equality follows from the multilinearity of $h$. Applying Lemmas~\ref{lem:bias} and~\ref{lem:analytic quotient} to this identity for $y\in G$ and using the hypothesis gives \eqref{eq:kernel-tail}.
\end{proof}
\subsection{Finding a tensor of small analytic rank}
Under the hypotheses of Proposition~\ref{prop:relativeschur}, we show how to obtain a form of small analytic rank by subtracting a form whose partition rank is bounded linearly in the number of defining equations.
\begin{proposition}\label{prop:all}
Let $Y$ and $Z$ be vector spaces, and let $T$ be a multilinear form on $X_1\times\cdots\times X_k\times Y\times Z$, where $k\geq1$. Suppose that $\Gamma_k$ consists of $s$ multilinear equations on the first $k$ blocks and that $\rk_{\F}M_T(x)\leq t$ on $V(\Gamma_k)$. There is a multilinear form $S$ on the same spaces with the following properties.
\begin{enumerate}
\item[$\mathrm{(a)}$] If $k=1$, then $\prank_{\F}(T-S)\leq s$ and $\arank_{\F}(S)\leq t$.
\item[$\mathrm{(b)}$] If $k\geq2$, then $\prank_{\F}(T-S)\ll_k s$ and
\begin{equation}\label{eq:all}
 \arank_{\F}(S)\ll_k \frac{s}{m}+m^{\,2^{k-1}-1}t
\end{equation}
for every integer $m$ satisfying $2(k-1)\leq m$ and $m^{2^{k-2}}\leq q$.
\end{enumerate}
\end{proposition}
\begin{proof}
Fix $k$ and put $C=k2^{k-1}$. Set $H_k=T$ and take $\Gamma_k$ as in
the hypothesis. We claim that for each $j=k,k-1,\ldots,1$, we can choose a system
$\Gamma_j$ of multilinear equations on $X_{[j]}$ and a multilinear
form $H_j$ on the original $k+2$ blocks with the following properties.
Write
\[
 s_j=|\Gamma_j|,\qquad
 b_j=\max_{x\in V(\Gamma_j)}\arank_{\F}\bigl(H_j(x,\cdot)\bigr).
\]
These quantities satisfy the following bounds.
\begin{enumerate}
\item[\textup{(i)}] $s_j+\prank_{\F}(T-H_j)\leq C^{k-j}s$.
\item[\textup{(ii)}] If $k\geq2$, then
$m^{2^{j-1}}b_j\leq C^{k-j}\bigl(s+m^{2^{k-1}}t\bigr)$ for every $m$ in
the range of \textup{(b)}.
\end{enumerate}

Assume the claim and take $j=1$.
The equations in $\Gamma_1$ are linear, so $L:=V(\Gamma_1)$
is a linear subspace of $X_1$ of codimension $c\leq s_1$. Choose a projection $p\colon X_1\to L$
and set $S(x_1,\cdot)=H_1(px_1,\cdot)$. If $(v_1,\ldots,v_c)$ is
a basis of $\ker p$, write $x_1-px_1=\sum_{i=1}^c\lambda_i(x_1)v_i$
with $\lambda_i\in X_1^*$. Then
\[
 (H_1-S)(x_1,\cdot)=\sum_{i=1}^c\lambda_i(x_1)H_1(v_i,\cdot).
\]
Each summand is a partition term, so by \textup{(i)}
$\prank_{\F}(T-S)\leq\prank_{\F}(T-H_1)+s_1\leq C^{k-1}s$.
For uniformly distributed $x_1\in X_1$, the vector $px_1$ is uniform on $L$. Thus
\[
q^{-\arank_{\F}(S)}
=\E_{u\in L}\,q^{-\arank_{\F}(H_1(u,\cdot))}
\geq q^{-b_1}.
\]
Hence $\arank_{\F}(S)\leq b_1$, which gives the desired bound in each case.

We prove the claim by descending induction on $j$. The case $j=k$ follows from the hypothesis.

Suppose that $(\Gamma_j,H_j)$ satisfies
\textup{(i)} and \textup{(ii)} for some $2\leq j\leq k$. Put $\ell=j-1$, write
$x=(x_1,\ldots,x_\ell)$ and $y=x_j$, and let $\Delta\subseteq\Gamma_j$
consist of the equations whose supports omit $y$. Arrange the remaining
$s_j-|\Delta|$ equations as the rows of a map
$B(x)\colon X_j\to\F^{s_j-|\Delta|}$. Each equation is linear in $y$ and multilinear on its support. Thus $B$ is affine in each parameter block, and
\[
 V(\Gamma_j)=\{(x,y):x\in V(\Delta),\ y\in\ker B(x)\},
 \qquad e:=\max_{x\in V(\Delta)}\rk_{\F} B(x)\leq s_j-|\Delta|.
\]
Apply Proposition~\ref{prop:kernel} to $H_j$, $\Delta$, and $B$ to obtain $H_{j-1}$ and $\Gamma_{j-1}$. The proposition gives
\begin{equation}\label{eq:rec-pr}
 s_{j-1}+\prank_{\F}(H_j-H_{j-1})
 \leq (2^\ell-1)|\Delta|+2^\ell e\leq 2^\ell s_j .
\end{equation}
By subadditivity, \eqref{eq:rec-pr} bounds
$s_{j-1}+\prank_{\F}(T-H_{j-1})$ by
$2^\ell s_j+\prank_{\F}(T-H_j)$. Since $2^\ell\leq C$, property
\textup{(i)} at $j$ gives \textup{(i)} at $j-1$.

Put $M=m^{2^{j-2}}$. Then $2\ell\leq M\leq q$ and
$M^2=m^{2^{j-1}}$. By \eqref{eq:kernel-tail},
\begin{align*}
 m^{2^{j-2}}b_{j-1}
 &\leq \frac{\ell 2^{\ell-1}M}{M-\ell}e
       +M\bigl(1+(2^\ell-1)(M-1)\bigr)b_j\\
 &\leq \ell 2^\ell e+2^\ell m^{2^{j-1}}b_j\\
 &\leq 2^\ell C^{k-j}\bigl((\ell+1)s+m^{2^{k-1}}t\bigr)\\
 &\leq C^{k-j+1}\bigl(s+m^{2^{k-1}}t\bigr),
\end{align*}
where we used $e\leq s_j\leq C^{k-j}s$ and $(\ell+1)2^\ell\leq C$.
This proves \textup{(ii)} at $j-1$ and completes the induction.
\end{proof}
\section{Putting everything together}\label{sec:main}
We now combine Propositions~\ref{prop:fibers}, \ref{prop:relativeschur}, and~\ref{prop:all} to prove Theorem~\ref{thm:partial}.
\begin{proof}[Proof of Theorem~\ref{thm:partial}]
If $T=0$, take $T'=0$. Suppose $T\ne0$. Write $k=d-2$ and $D=2^{k-1}$, and for now assume $q\geq(2k)^D$. We view the first $k$ blocks as $X=X_1\times\cdots\times X_k$ and denote the last two by $Z,W$.

Proposition~\ref{prop:fibers} gives a system $\Gamma$ of multilinear equations on $X$ such that
\[
 s_0=|\Gamma|,\qquad b=\max_{x\in V(\Gamma)}\rk_{\F}M_T(x),\qquad s_0+b\leq K_{k+1}(r+1).
\]
Since $q\geq 2k>k$, Proposition~\ref{prop:relativeschur} gives a form $S$ and a refined system $\Gamma'$ of $s$ multilinear equations on $X$ such that
\[
 \prank_{\F}(T-S)\leq 2^k b,\qquad s\leq(2^k-1)s_0,\qquad
 \max_{x\in V(\Gamma')}\rk_{\F}M_S(x)\leq t:=\frac{k2^{k-1}b}{q-k}.
\]
The inequality $q-k\geq q/2$ gives $t\leq k2^k b/q$.

For $k=1$, Proposition~\ref{prop:all}(a) applied to $S$ and $\Gamma'$ gives a form $T'$ such that $\prank_{\F}(S-T')\leq s$ and $\arank_{\F}(T')\leq t\leq 2b/q$. Since $D=1$, this gives the required power of $q$.

For $k\geq2$, Proposition~\ref{prop:all}(b) gives a form $T'$ with $\arank_{\F}(T')\ll_k s/m+m^{D-1}t$. To bound both terms by a multiple of $(s+b)q^{-1/D}$, take $m=\lfloor q^{1/D}\rfloor$. Indeed, $q^{1/D}\geq2k$, so $m\geq2(k-1)$. Also $m^{2^{k-2}}=m^{D/2}\leq q^{1/2}\leq q$. This verifies both restrictions on $m$. Since $q^{1/D}\geq2$, we have $m\geq q^{1/D}/2$, whereas $m^{D-1}\leq q^{(D-1)/D}$. Therefore
\[
 \frac{s}{m}\leq 2s q^{-1/D},\qquad
 m^{D-1}t\leq q^{(D-1)/D}\frac{k2^k b}{q}=k2^k b q^{-1/D}.
\]
It follows that $\prank_{\F}(S-T')\ll_k s$ and $\arank_{\F}(T')\ll_k(s+b)q^{-1/D}$, and these bounds also hold for $k=1$.

By subadditivity, the partition rank of $T-T'$ is at most the sum of the bounds for $T-S$ and $S-T'$, which is $O_k(b+s)$. Moreover, $b+s\leq (2^k-1)(b+s_0)\leq (2^k-1)K_{k+1}(r+1)$. Since $k=d-2$, we may choose constants $L_d^0,B_d\geq1$ depending only on $d$ such that
\[
 \prank_{\F}(T-T')\leq L_d^0(r+1),\qquad
 \arank_{\F}(T')\leq B_d(r+1)q^{-1/D}.
\]

Set $L_d=3L_d^0$ and choose an integer $Q_d$ at least $\max\{(2k)^D,(d-1)^2,(6B_d)^D\}$. For $q\geq Q_d$, the first threshold ensures that the preceding construction applies. The second gives $r\geq1/2$ by Lemma~\ref{lem:bias}, since $T\ne0$. Thus $r+1\leq3r$, and our partition-rank bound is at most $L_dr$. Finally, the third threshold gives $3B_d q^{-1/D}\leq1/2$, so the analytic rank of the same $T'$ is at most $B_d(r+1)q^{-1/D}\leq r/2$. This proves \eqref{eq:partial}.
\end{proof}
\subsection*{Declaration of AI use}
Multiple versions of ChatGPT (mainly 5.6 Sol and 6 Astra) were used by the authors to develop the idea of the proof. The AI tool first found a proof of the conjecture for 4-tensors(especially the 4-tensor version of Proposition~\ref{prop:all}). The authors then generalized the proof to any order tensors. All AI-assisted outputs were independently checked, verified, and simplified by the authors, who assume full responsibility for the content of this work.
\bibliographystyle{abbrv}
\bibliography{ref}
\end{document}